\documentclass[10pt,a4paper,final]{article}
\usepackage[margin=0.5in]{geometry}
\usepackage{amsmath}
\usepackage{amsfonts}
\usepackage{amssymb}
\usepackage{hyperref}
\usepackage[final]{graphicx}

\renewcommand{\vec}[1]{\ensuremath{\mathbf{#1}}}

\newcommand{\mnm}{\ensuremath{\text{minimize}}}
\newcommand{\wrt}{\ensuremath{\mathrm{with \; respect \; to}}}
\newcommand{\sto}{\ensuremath{\mathrm{subject \; to}}}

\title{Sizing Optimization of Truss Structures using a Hybridized Genetic Algorithm}
\begin{document}
 \begin{titlepage}
\maketitle
{
\centering
Reza Najian Asl{\bf *}, \footnote[1]{{\bf *}Corresponding author: Research Assistant, Mechanical Engineering Department, University of Tehran, Tehran, Iran. \href{reza.najian-asl@tum.de}{reza.najian-asl@tum.de}}, Research Assistant\\
Mohamad Aslani, \footnote[2]{Research Assistant, Mechanical Engineering Department, University of Tehran, Tehran, Iran. \href{mailto:maslani@iastate.edu}{maslani@iastate.edu}}, Research Assistant \\
Masoud Shariat Panahi \footnote[3]{Mechanical Engineering Department, University of Tehran, Tehran, Iran. \href{mailto:mshariatp@ut.ac.ir}{mshariatp@ut.ac.ir}}, Associate Professor\\
}
\end{titlepage}
\abstract
This paper presents a genetic-based hybrid algorithm that combines the exploration power of Genetic Algorithm (GA) with the exploitation capacity of a phenotypical probabilistic local search algorithm. Though not limited to a certain class of optimization problems, the proposed algorithm has been ``fine tuned" to work particularly efficiently on the optimal design of planar and space structures, a class of problems characterized by the large number of design variables and constraints, high degree of non-linearity and multitude of local minima. The proposed algorithm has been applied to the skeletal weight reduction of various planar and spatial trusses and shown to be superior in all of the cases.

\section{Introduction}
Skeletal structures are widely used as benchmark test in structural optimization as they can include a large number of design variables and constraints leading to high degree of non linearity and presence of local minima. The demand for reliable, computationally inexpensive optimum structural design tools has motivated the researchers to develop specialized optimization techniques and/or to ``tune" existing methods to solve this class of problems more efficiently.\\
Structural design problems are typically characterized by their large numbers of design variables and constraints, high degree of non-linearity and multitude of local minima. It is therefore extremely difficult, if not practically impossible, to find the globally optimum solution to the problem of designing large scale space structures unless some additional knowledge of the optimum's whereabouts is available to guide the search. If this specific knowledge is not available, one has to resort to search techniques that can, to a certain extent, escape the local minima and spot the global one(s).\\
It comes with no surprise that structural optimization has grown to become a challenging one, seeking to determine the structure's dimensions, geometry and topology that would render the structure as
light/inexpensive as possible while keeping its performance characteristics (e.g. stresses and displacements) within allowable limits. The inherent complexity of the problem has rendered it a perfect benchmark for large scale search algorithms.\\
It is widely believed that Evolutionary Algorithms (EAs) are so far the most promising methods for searching large, highly non-linear design space black box problems with many local minima. EAs are stochastic search
methods that mimic the metaphor of natural biological evolution. They operate on a population of potential solutions applying the principle of survival of the fittest to produce increasingly better
approximations to a solution. Through the adaptation of successive generations of a large number of individuals, an evolutionary algorithm performs an efficient directed search. Evolutionary search is generally better than random search and is not susceptible to the hill-climbing behaviours of gradient-based search \cite{Sivanandam2008}.\\
The most common EA is Genetic Algorithm (GA) that was introduced by Holland in 1975 \cite{Holland1992}. In Holland's view, GAs are a computational analogy of adaptive systems. They are modelled loosely on the principles of the evolution via natural selection, employing a population of individuals that undergo selection in the presence of variation-inducing operators such as mutation and recombination(crossover).
A fitness function is used to evaluate individuals, and reproductive success varies with fitness \cite{Sivanandam2008}.
From a practical point of view, traditional GAs suffer from two problems: their high computational cost and their vulnerability to getting trapped in local minima of the objective function, especially when the number of variables is large and/or the design space is bounded with a huge number of constraints. Many modifications have been performed on the operators of GAs to resolve these problems, most of them
aimed at changing GAs to increase their rate of convergence. A more recent trend has been to hybridize GAs with other, more specialized, optimization techniques such as Simulated Annealing (SA) \cite{Bolte1996} and \cite{Wang2004}, Nelder-Mead(Simplex) \cite{Wren2007}, \cite{Rahami2011} and \cite{Aslani2010} and Monte Carlo algorithm \cite{N.Dugan2009} and particle swarm optimizer \cite{Fan2006}.\\
Since the standard GA operators are traditionally unable to secure a local search in the neighbourhood of existing solutions, it might be argued that hybridization of GA with a neighbourhood search
algorithm would improve its convergence rate. However, this neighbourhood search algorithm must not undermine the stochastic nature of GA and must be justifiably simple so that that the resulting hybrid
algorithm would not border with a simple random walk in the design space. These conditions make SA a prime candidate for the job; a simple, yet powerful stochastic search algorithm that is not fooled by false minima and is easy to implement.
However, Simulated Annealing, even after recent improvements aimed at boosting its convergence behaviour \cite{Carbas2010}, comes with its own drawbacks, including its high sensitivity to the location of initial points. The farther from the global optimum they are, the more iteration it takes the algorithm to get there. This is because the static neighbourhood search of the SA(one with a constant neighbourhood radius throughout the search) makes the algorithm unable to adjust its neighbourhood size according to its distance from the optimum. This could be avoided by introducing a ``dynamic" neighbourhood search where the radius of the local search area is adaptively varied according to some measure of ``closeness" to the global optimum.
In the next section, a new hybrid GA is introduced which adopts a modified version of SA's stochastic neighbourhood search to improve its convergence rate. Then back to our structural sizing optimization problem and it is shown that how the problem could be solved using the proposed algorithm.
\section{The Hybrid Simulated Annealing-Genetic Algorithms(H-SAGA)}
In order to clarify nomenclature, the basis of GA and SA algorithms included in the present optimization code will be briefly recalled.
\subsection{The Basic Genetic Algorithm}
In a Genetic Algorithm, a predetermined number(a population) of strings(chromosomes) which encode candidate solutions to an optimization problem evolves toward better solutions. Actual solutions(phenotype) are represented by (binary/real value) encoded strings(genotype). The evolution usually starts from a population of randomly generated individuals and happens in generations. In each generation, the fitness(measure of desirability) of every individual in the population is evaluated.\\
Multiple individuals are then stochastically selected from the current population (based on their fitness), and modified using such genetic operators as recombination(cross over) and mutation to form a new population. The new population is then used as the ``current population" in the next iteration of the algorithm.\\
The population size depends on the nature of the problem. While a large population size tends to secure the diversity of the candidate solutions and avoids a premature convergence, it inevitably increases
the computational cost of the search. Special arrangements are usually made to encourage a relatively uniform distribution of individuals across the search space. Occasionally, the solutions may be ``seeded" in areas where optimal solutions are likely to be found.\\
To breed a new generation, individual solutions are selected to form a ``mating pool" through a fitness based process, where fitter solutions (as measured by a fitness function) would have a higher chance to be
selected. In GA the selection mechanism is always partially (if not totally) stochastic to allow some proportion of less fit solutions be possibly selected. This helps keep the diversity of the population large, preventing premature convergence to poor solutions.\\
For each new solution to be produced, a pair of ``parent" solutions is selected for breeding from the mating pool. By producing a ``child" solution using genetic operators, a new solution is created which
typically shares many of the characteristics of its ``parents". New parents are selected for each new child, and the process continues until a new population of solutions of appropriate size is generated. The average fitness of the new generation is expected to be higher than that of the previous generation, as the best individuals of the last generation have been given higher chances for breeding.\\
The algorithm terminates when either a maximum number of generations has been produced, or a satisfactory fitness level has been reached for the population.
\subsection{The Basic Simulated Annealing}
Initially developed in the field of statistical mechanics to simulate a collection of atoms in equilibrium at a given temperature\cite{Kirkpatrick1983}, Simulated Annealing was soon recognized as a powerful optimization technique. The name and inspiration come from the ``annealing" process in metallurgy, a technique involving heating and controlled cooling of a material to increase the size of its crystals and reduce their defects. The heat causes the atoms to become unstuck from their initial positions(a local minimum of the internal energy) and wander randomly through states of higher energy; the slow cooling gives them more chances of finding configurations with lower internal energy than the initial one.
By analogy with this physical process, each step of the SA algorithm replaces the current solution($A$) by a random ``nearby" solution($B$), generated within a neighbourhood of predetermined size, with a probability($P$) that depends on the difference between the corresponding function values and on a global parameter $T$ (called the temperature), that is gradually decreased during the process. This probability, in its basic form, is calculated as follows:\\
\begin{equation}
	\begin{aligned}
		\label{Prob}
  		P = \left \{ 
  		\begin{array}{l l}
    		1(i.e\quad 100 \%) & \quad \text{if f(A) $\geq$ f(B)}\\
    		exp(\frac{f(A)-f(B)}{T} ) & \quad \text{otherwise}\\
  		\end{array} \right.
  	\end{aligned}
\end{equation}
This formula(Eq.\ref{Prob}) is commonly justified by analogy with the transitions of a physical system. The essential property of this probability function is that it is non zero when $f(B) \geq f(A)$, meaning that the system may move to the new state even when it is worse (has a higher energy) than the current one. It is this feature that prevents the algorithm from becoming stuck in a local minimum (a state that is worse than the global minimum, yet better than any of its neighbours).\\
The algorithm terminates when either the temperature falls below a certain value or the objective function does not improve much during a certain number of consecutive iterations. As mentioned earlier, an improved SA formulation will be used in which the neighbourhood of local search is adaptively re-sized throughout the search.\\
\subsection{Handling Constraints}
A constrained problem in which the feasible region is defined by a set of implicit/explicit constraints could be treated in two ways. One is to use them as preventive measures that would not allow a newly generated solution to enter the population unless it satisfies all the constraints; and the other is to use them as penalizers, that is to allow infeasible solutions to enter the population but penalize (increase, in case of a minimization problem) their fitness in proportions with their degree of constraint violation. The latter is usually preferred because it provides the user with a quantitative measure of the infeasibility of a solution and helps the algorithm find its way to the feasible region.\\
In structural optimization, weight of the structure usually constitutes the objective function and cross-sectional properties of the members form the variable set. Also, the structure is constrained in terms of the stresses and displacements of its members. These are considered indirect constraints as they are not applied directly on the cross sections of the members.\\
These non-linear constraints are incorporated into the objective function to form a new, unconstrained function. In this way, infeasible solutions are given a varying amount of penalty to subside their fitness. In an attempt to make the search somewhat adaptive, dynamic penalties are introduced that would vary with the number of iterations. The problem could then be formulated as follows (without loss of generality, only the inequality constraints are considered here):
\begin{equation}
	\begin{aligned}
		\label{Prob}
  		\mnm &\quad f(\vec{x})		\\
  		\wrt &\quad \vec{x} = (x_1,\cdots, x_n) \\
		\sto &\quad g_i(\vec{x}) \leq 0 , \forall \quad i = 1 \quad \text{to} \quad m \\
  		 &\quad					\\
  	\end{aligned}
\end{equation}
The new, unconstrained objective function would then be defined as:\\
\begin{equation}
	\begin{aligned}
		\label{ProbP}
		&\quad F(\vec{x}) = f(\vec{x})\pm p(\vec{x})
	\end{aligned}
\end{equation}
where\\
\begin{equation}
	\begin{aligned}
		\label{ProbP1}
		p(\vec{x})= \alpha (\# iteration)^\beta \times \sum_{i = 0}^q S_i(\vec{x})
	\end{aligned}
\end{equation}
and
	\begin{equation}
		\begin{aligned}
			\label{ProbP2}
			S_i(\vec{x}) = \left\{						
			\begin{array}{l l}
			0 &\quad g_i(\vec{x}) \leq 0 \\
			\| g_i(\vec{x}) \| &\quad otherwise
			\end{array} \right.
	\end{aligned}
\end{equation}

\subsection{The Dynamic Neighbourhood Search (DNS)\label{DNS}}
To accelerate the search and reduce the number of times the algorithm falls into and climbs out of local valleys (in a minimization problem), the neighbourhood of local search for each design variable is initially set equal to the difference between its highest and lowest values within the 10 best individuals, a large value meant to encourage exploratory leaps. The neighbourhood is then progressively reduced using an improvement-dependent criterion aimed at refining the search steps once large exploratory steps could no longer improve the solution. The performance of DNS is controlled by two user-defined parameters: the radius reduction threshold $\beta$ and the contraction parameter $\gamma$. The implementation of DNS is somewhat elaborated in the pseudo code below.

\subsection{The Hybrid Algorithm \label{THA}}
In essence, the proposed hybrid algorithm works as follows:
\begin{itemize}
\item	An initial population of encoded potential solutions is randomly generated (this may contain infeasible solutions);
\item	Fitness of each solution is calculated using the alternative objective function (Eq. \ref{ProbP});
\item	The population is evolved according to GA;
\item	Every $T_{SA}$ (Simulated Annealing Threshold) generations, a SA search is performed with its initial point set at the fittest individual in generation $T_{SA}$. The search region is defined by current positions of the 10 fittest solutions in generation $T_{SA}$ and the size of the neighbourhood to be locally searched is calculated according to section \ref{DNS}.  The SA search is terminated once the average fitness improvement within a certain number of consecutive iterations falls below a search precision ($\epsilon$) and or the temperature falls below a predefined value ($T_{MIN}$) implying that further probabilistic moves to inferior neighbouring points would be highly unlikely;
\item	For every new solution suggested by SA, an individual is removed from the population with a probability proportional to the inverse of its fitness;
\item	The algorithm terminates once the termination criteria for the GAs is met;
\end{itemize}
A pseudo code for the algorithm is presented below:\\
\begin{center}
\line(1,0){250}
\end{center}
{\bf Choose} GA parameters, radius reduction threshold (TSA)\\
{\bf Generate} initial population (randomly)\\
{\bf Do}\\
\hspace*{10mm}	{\bf Evaluate} each individual's constraint violation (Eq. \ref{ProbP1})\\
\hspace*{10mm}	{\bf Evaluate} unconstrained objective function (Eq. \ref{ProbP}) \\
{\bf Form} the new generation\\
\hspace*{10mm}	{\bf Use} fitness-proportionate selection to form the mating pool\\
\hspace*{10mm}	{\bf Mate} pairs at random\\
             \hspace*{20mm} {\bf Apply} crossover operator\\
             \hspace*{20mm} {\bf Apply} mutation operator\\
             \hspace*{20mm} {\bf Rank } individuals based on their fitness\\
             \hspace*{20mm} {\bf Select} 10 best-ranking individuals\\
             \hspace*{25mm} {\bf If} number of generations since last application of SA equals $T_{SA}$\\
             \hspace*{25mm} {\bf Apply} SA\\
             \hspace*{29mm} {\bf Define} an initial temperature $T$\\
             \hspace*{29mm} {\bf Set} the search radius of each design variable ($R_i$) (Section \ref{DNS})\\
			 \hspace*{29mm} {\bf Choose} the cooling coefficient $\alpha$ ($ 0< \alpha < 1$)\\
			 \hspace*{29mm} {\bf Choose} radius reduction threshold $\beta$ and the contraction parameter $\gamma$\\
			 \hspace*{29mm} {\bf Choose} the top-ranking individual as the starting point\\
			 \hspace*{29mm} {\bf Do}\\
 \hspace*{29mm} {\bf Choose} a random point ($N$) in the neighbourhood of current point($C$)\\
 \hspace*{29mm} {\bf Set} $\delta = f(P)-f(C)$\\
 \hspace*{29mm} {\bf If} $\delta < 0 $ by probability $P$ (Eq. \ref{Prob})\\
 \hspace*{29mm} {\bf Replace} $C$ with $N$\\
 \hspace*{29mm} {\bf If} number of SA iterations equals number of design variables  $T = \alpha T$\\
 \hspace*{29mm} {\bf Endif}\\
 \hspace*{29mm} {\bf If} iterations falls below $\beta$\\
 \hspace*{20mm} {\bf Then} upgrade search radius $R_{(i,new)} = R_{(i,old)}/ \gamma$\\
 \hspace*{25mm} {\bf Endif} \\
 \hspace*{20mm} {\bf While} SA termination criteria met (Section \ref{THA})\\
 \hspace*{20mm} {\bf Add} SA solution to the population and remove an individual based on the inverse of their fitness\\ 
 \hspace*{25mm} {\bf Endif}\\
 {\bf While} GA termination criteria is met \\
\begin{center}
\line(1,0){250}
\end{center}

\section{Application to Structural Design Problems}
A variety of structural design problems ranging from a 10-bar planar truss to a 26-story space structure with various numbers of variables and constraints are considered in this section.\\
From an analytical point of view, the problem of weight minimization of a truss structure is defined as:\\
\begin{equation}
	\begin{aligned}
		\label{ProbMain}
  		\mnm &\quad W(\vec{A}) = \rho g \sum_{j = 1}^{NEL} A_j l_j		\\
  		\wrt &\quad \vec{A} = (A_1,\cdots, A_{NEL}) \\
		\sto &\quad \left\{ 
  		\begin{array}{l l}
  		u_{(x,y,z),k}^l \leq u_{(x,y,z),k, ilc} \leq u_{(x,y,z),k}^u\\
  		\sigma_j^l \leq \sigma_{j,ilc} \leq \sigma_j^u &\quad  \left\{ 
  		\begin{array}{l l}
		i = 1, NEL\\
		k = 1, NOC\\
		l = 1, NLC\\
  		\end{array} \right.\\
  		A_j^l \leq A_j \leq A_j^u\\
  		\end{array} \right.
  	\end{aligned}
\end{equation}
Where:
\begin{itemize}

\item	$A_i$ are the cross sectional areas of the members;
\item	$W$ is the total weight of the structure;
\item	$NEL$ is the number of elements in the structure;
\item	$NOD$ is the total number of nodes;
\item	$NLC$ is the number of independent loading conditions acting on the structure;
\item	$g$ is the gravity acceleration;
\item	$\rho$ is the material density used in the structure;
\item	$\vec{A}$ is the vector of the cross sectional area values; $A_j$ is the \emph{$j_{th}$} element area of the structure.\item		$A_j^l$  and $A_j^u$ are lower and upper bounds of $A_j$;
\item	$l_j$ is the length of \emph{$jth$} element, and is calculated according to Eq. \ref{length};
\begin{equation}
\label{length}
	\begin{aligned}
		l_j= \sqrt{(x_{j1}-x_{j2} )^2+(y_{j1}-y_{j2} )^2+(z_{j1}-z_{j2})^2} 
	\end{aligned}		
\end{equation}
	
where $x_{j1,2},y_{j1,2},z_{j1,2}$ are coordinates of the nodes of the $jth$ element;
\item $u_{(x,y,z),k,ilc}$ is  the displacement of the $kth$ node in the directions $x,y,z$ with the lower and upper limits $u_{(x,y,z),k}^l$ and $u_{(x,y,z),k}^u$ in corresponding of the $ilcth$ loading condition, Eq. \ref{ProbMain};
\item $\sigma_{j,ilc}$ is the stress on the $jth$ element, $\sigma_k^l$ and $\sigma_k^u$ are the lower and upper limits in correspondence of the $ilcth$ loading condition, Eq.\ref{ProbMain};
\end{itemize}

Design variables ($A_i$) are treated as continuous variables, as it is a common practice in the literature. Each truss structure is analysed using the \emph{FEM displacement method} and the results are presented in a separate table. Attention must be paid to the fact that the numbers reported in these tables have been rounded up to 4 digits. They may therefore seem to slightly violate some of the constraints or deviate from global optimality, a problem that does not arise when the actual, 7-digit values are considered.

\subsection{The 10-bar truss}
\begin{figure}[thb]
\begin{center}
  \centering
 \includegraphics[scale=1]{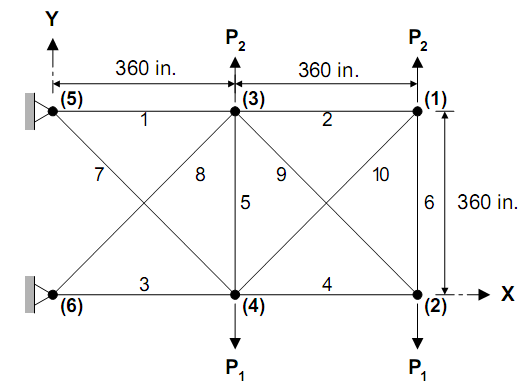} 
\caption{10-bar planar truss}
\label{fig:10bar}
\end{center}
\end{figure}

The cantilever truss shown in Fig. \ref{fig:10bar} with 10 independent design variables has been studied by many researchers Schmit and Farshi \cite{L.A.Schmit1974}, Schmit and Miura \cite{L.A.Schmit1976},Venkayya \cite{Venkayya1971}, Sedaghati \cite{Sedaghati2005}, Kaveh and Rahami \cite{Kaveh2006}, Li et al. \cite{Li2007}, Farshi and Ziazi \cite{Li2007}, Rizzi \cite{RizziMay1976}, John et al. \cite{John1987}. A material density of $0.1 lb/(in^3 )$ and a modulus of elasticity of 10000 $ksi$ is commonly used. Members are constrained in displacement and stress to $\pm2 in$ and $\pm25 ksi$ respectively. The lower bound of cross sectional areas is $0.1 in^2$. Two load cases are considered (Fig.\ref{fig:10bar}):
\begin{enumerate}
\item A single load ($P_1=100 kips$ and $P_2=0 kips$) 
\item Double loads ($P_1=150 kips$ and $P_2=50 kips$)
\end{enumerate}

Tables 1 and 2 give the optimum design for Cases 1 and 2, respectively, and also provide a comparison between the optimal design results reported in the literature and the present work\footnote{The results included here and in the rest of this paper for comparison purposes are the ``feasible" ones. A few reported solutions that involve constraint violation(s), including those reported in references \cite{Li2007} and \cite{Lee2004} are not included for some cases here.}.\\
\begin{table}
\caption{Optimal design comparison for the 10-bar planar truss (Case 1)}
\centering
\resizebox{\textwidth}{!}{%
\begin{tabular}{l*{7}{c}r}
Variables($in^2$)   & Schmit and Farsh \cite{L.A.Schmit1974} & Schmit and Miura \cite{L.A.Schmit1976} & Venkayya \cite{Venkayya1971} & Sedaghati \cite{Sedaghati2005} & Kaveh and
Rahami \cite{Kaveh2006}  & Li et al. \cite{Li2007} & Farshi and Ziazi \cite{Farshi2010} & This work  \\ \hline
$A_1$ 			& 33.43 & 30.67 & 30.42 & 30.5218 & 30.6677 & 30.704 & 30.5208 & 30.5091 \\
$A_2$                                & 0.100 & 0.100 & 0.128 & 0.1000 & 0.1 &  0.100 & 0.1000 &  0.1000  \\
$A_3$                                & 24.26 & 23.76 & 23.41 & 23.1999 & 22.8722 &  23.167 & 23.2040 &  23.2004  \\
$A_4$     		          & 14.26 & 14.59 & 14.91 & 15.2229 & 15.3445 & 15.183 & 15.2232 & 15.1926  \\
$A_5$			          & 0.100 & 0.100 & 0.101 & 0.1000 & 0.1 & 0.100 & 0.1000 & 0.1000\\
$A_6$			          & 0.100 & 0.100 & 0.101 & 0.5514 & 0.4635 & 0.551 & 0.5515 & 0.5559\\
$A_7$			          & 8.388 & 8.578 & 8.696 & 7.4572 & 7.4796 & 7.460 & 7.4669 & 7.4612\\
$A_8$			          & 20.74 & 21.07 & 21.08 & 21.0364 & 20.9651 & 20.978 & 21.0342 & 21.0714\\
$A_9$			         & 19.69 & 20.96 & 21.08 & 21.5284 &  21.7026 & 21.508 & 21.5294 & 21.4731\\
$A_{10}$		         & 0.100 & 0.100 & 0.186 & 0.1000 & 0.1 & 0.100 & 0.1000 & 0.1000\\ \hline
Weight$lb$		         & 5089.0 & 5076.85 & 5084.9 & 5060.85 & 5061.90 & 5060.92 & 5061.4 & 5058.66\\ \hline
Note: 1 $in^2 = 6.425 cm^2$ $1 lb = 4.448 N$
\label{T10bar2}
\end{tabular}}
\end{table}
\begin{table}
\caption{Optimal design comparison for the 10-bar planar truss (Case 2)}
\centering
\resizebox{\textwidth}{!}{%
\begin{tabular}{l*{7}{c}r}
Variables($in^2$)   & Schmit and Farsh \cite{L.A.Schmit1974} & Schmit and Miura \cite{L.A.Schmit1976} & Venkayya \cite{Venkayya1971} & Rizzi \cite{RizziMay1976} & John et al. \cite{John1987}  & Lie at al. \cite{Li2007} & Farshi and Ziazi \cite{Farshi2010} & This work  \\ \hline
$A_1$ 			& 24.29 & 23.55 & 25.19 & 23.53 & 23.59 & 23.353 & 23.5270 & 23.3187\\
$A_2$                                & 0.100 & 0.100 & 0.363 & 0.100 & 0.10 & 0.100 & 0.1000 & 0.1\\
$A_3$                                & 23.35 & 25.29 & 25.42 & 25.29 & 25.25 & 25.502 & 25.2941 & 25.5790\\
$A_4$     		          &13.66 & 14.36 & 14.33 & 14.37 & 14.37 & 14.250 & 14.3760 & 14.6640\\
$A_5$			          & 0.100 & 0.100 & 0.417 & 0.100 & 0.10 & 0.100 & 0.1000 & 0.1\\
$A_6$			          &1.969 & 1.970 & 3.144 & 1.970 & 1.97 & 1.972 & 1.9698 & 1.9695\\
$A_7$			         &12.67	& 12.39 & 12.08 & 12.39 & 12.39 & 12.363 & 12.4041 & 12.2654\\
$A_8$			         &12.54	 & 12.81 & 14.61 & 12.83 & 12.80 & 12.894 & 12.8245 & 12.6473\\
$A_9$			         &21.97 & 20.34 & 20.26 & 20.33 & 20.37 & 20.356 & 20.3304 & 20.3422\\
$A_{10}$		        & 0.100 & 0.100 & 0.513 &	0.100 & 0.10 & 0.101 & 0.1000 & 0.1\\ \hline
Weight$lb$		       & 4691.84 & 4676.96 & 4895.60 & 4676.92 & 4676.93 & 4677.29 & 4677.8 & 4675.43\\ \hline
Note: 1 $in^2 = 6.425 cm^2$ $1 lb = 4.448 N$
\label{T10bar2}
\end{tabular}}
\end{table}

\subsection{The 17-bar truss}
\begin{figure}[thb]
\begin{center}
  \centering
  \includegraphics[scale=1]{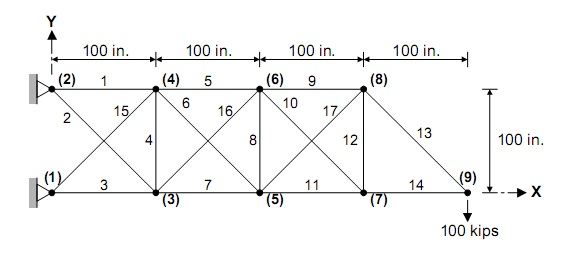} 
\caption{17-bar planar truss}
\label{fig:17bar}
\end{center}
\end{figure}

The 17-bar truss of Fig. \ref{fig:17bar} is also a typical test case studied by many researchers, including Khot and Berke \cite{Venkayya1973}, Adeli and Kumar \cite{H.Adeli1995}, Lee and Geem \cite{Lee2004}. In this case the material density is $0.268 lb/in^3$  and the modulus of elasticity is $30000 ksi$. A displacement limitation of $\pm2 in$ were the only constraint that imposed on the nodes in both $x$ and $y$ directions. A single loading of $100 kips$ is on the node $(9)$, Seventeen independent variables exist as there is no variable grouping in the problem. A minimum cross-section of $0.1 in^2$  is the lower boundary of the variables.\\
Table \ref{T17bar} compares the optimization results with similar studies published in literature. It can be seen that H-SAGA designed a lighter structure.
\begin{table}
\caption{Optimal design comparison for the 17-bar planar truss}
\centering
\resizebox{\textwidth}{!}{%
\begin{tabular}{l*{5}{c}r}
Variables($in^2$)   & Khot and Berke  \cite{Venkayya1973} & Adeli and Kumar \cite{H.Adeli1995} & Lee and Geem  \cite{Lee2004} & This work  \\ \hline
$A_1$ 			 & 15.93 & 16.029 & 15.821 & 15.8187 \\
$A_2$                                 & 0.1 & 0.107 & 0.108 & 0.1051\\
$A_3$                                & 12.07 & 12.183 & 11.996 & 12.0246\\
$A_4$     		          & 0.1 & 0.1 & 0.1 &  0.1\\
$A_5$			          & 8.067 & 8.417 & 8.15 & 8.1132\\
$A_6$			         & 5.562 & 5.715 & 5.507 & 5.5318\\
$A_7$			         & 11.933 & 11.331 & 11.829 & 11.8431\\
$A_8$			         & 0.1 & 0.105 & 0.1 & 0.1\\
$A_9$			        & 7.945 & 7.301 & 7.934 & 7.9560\\
$A_{10}$		        & 0.1 & 0.115 & 0.1 & 0.1\\ 
$A_{11}$                        & 4.055 & 4.046 & 4.093 & 4.0711\\
$A_{12}$                        & 0.1 & 0.101 & 0.1 & 0.1\\
$A_{13}$                        & 5.657 & 5.611 & 5.66 & 5.6841\\
$A_{14}$                       & 4 & 4.046 & 4.061 & 4.0087\\
$A_{15}$                      & 5.558 & 5.152 & 5.656 & 5.5849\\
$A_{16}$                      & 0.1 & 0.107 & 0.1 & 0.1\\
$A_{17}$                      & 5.579 & 5.286 & 5.582 & 5.5804\\ \hline
Weight$lb$		     & 2581.89 & 2594.42 & 2580.81 & 2578.76\\ \hline
Note: 1 $in^2 = 6.425 cm^2$ $1 lb = 4.448 N$
\label{T17bar}
\end{tabular}}
\end{table}
\subsection{The 18-bar truss}
\begin{figure}[thb]
\begin{center}
  \centering
  \includegraphics[scale=1]{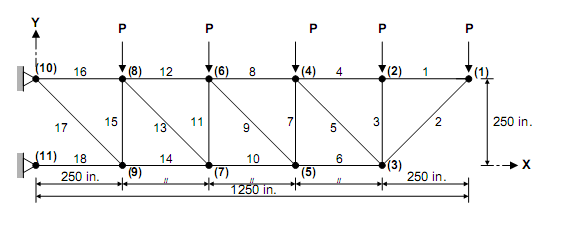} 
\caption{18-bar planar truss}
\label{fig:18bar}
\end{center}
\end{figure}

The optimum design of the 18-bar cantilever truss shown in Fig.\ref{fig:18bar} is a sizing optimization problem already analyzed by Imai and Schmit \cite{Imai1981} and Lee and Geem \cite{Lee2004}. Material density is $0.1 lb/(in^3 )$ . Stress in members must not exceed $20000 psi$ (The same limit in tension and compression). Constraints on Euler buckling strength are also considered where the allowable limit for the $ith$ member is computed as: 
\begin{equation}
\label{length}
	\begin{aligned}
	\sigma^b_i = - \frac {KEA_i}{L_i^2}
	\end{aligned}
\end{equation}
where $K$ is the buckling constant ($K=4$). $E$ is the modulus of elasticity($E=10000 ksi$) and  $L_i$ is the length of the element. Downward load of $F=20 kips$ is on nodes (1), (2), (4), (6) and (8). The member cross-section groups are as follow:
\begin{enumerate}
\item $A_1 = A_8 = A_{12} = A_4 = A_{16}$
\item $A_2 = A_6 = A_{10} = A_{14} = A_{18}$  
\item $A_3 = A_7 = A_{11} = A_{15}$
\item $A_5 = A_9 = A_{13} = A_{17}$
\end{enumerate}
Minimum cross section is $0.1 in^2$. H-SAGA obtained a lighter designing in comparison with other algorithms. These results are presented in Table \ref{T18bar2}.
\begin{table}
\caption{Optimal design comparison for the 18-bar planar truss}
\centering
\begin{tabular}{l*{4}{c}r}
Variables($in^2$)   & Imai and Schmit  \cite{Imai1981} & Lee and Geem  \cite{Lee2004} & This work  \\ \hline
$A_1$ 			& 9.998 & 9.98 & 9.9671\\
$A_2$                                & 21.65 & 21.63 & 21.5990\\
$A_3$                                & 12.5	 & 12.49 & 12.4492\\
$A_4$     		          & 7.072 & 7.057 & 7.0490\\ \hline
Weight$lb$		         & 6430 & 6421.88 & 6419.23\\ \hline
\label{T18bar2}
\end{tabular}
\end{table}
\subsection{The 22-bar space truss}
\begin{figure}[thb]
\begin{center}
  \centering
  \includegraphics[clip, width = 120mm , height = 75mm]{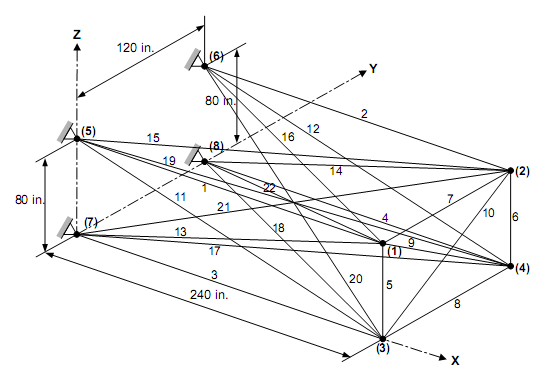} 
\caption{22-bar space truss}
\label{fig:22bar}
\end{center}
\end{figure}

Weight minimization of the spatial 22-bar truss shown in Fig.\ref{fig:22bar} was previously attempted with or without layout variables by Sheu and Schmit \cite{Sheu1972}, Khan and Willmert \cite{Khan1979}, Li et al \cite{Li2007}, Farshi and Ziazi \cite{Farshi2010}. Material density was taken as $0.1 lb/in^3 $ , the modulus of elasticity as $10000 psi$. Members were grouped as follows:\\
$(1)A_1 \sim A_4 (2) A_5 \sim A_6 (3) A_7 \sim A_8(4)A_9\sim A_{10}(5)A_{11} \sim A_{14} (6) A_{15} \sim A_{18}(7) A_{19} \sim A_{22}$
Stress members are limited according to Table \ref{T22bar1}. Nodal displacement must be less than $2 in$. The structure is subject to three independent loading conditions (Table \ref{T22bar2}). The lower bound of cross sectional areas is $0.1 in^2$. Optimization results are presented in Table \ref{T22bar3}.
\begin{table}
\caption{Stress limits for the spatial 22-bar truss problem}
\centering
\begin{tabular}{l*{4}{c}r}
Variables($in^2$)   & Grouping & Compressive Stress Limitation (ksi) & Tensile stress limitation(ksi) \\ \hline
1 & $A_1\sim A_4$ 		 & 24.0 & 36.0  \\
2 & $A_5 \sim A_6$          & 30.0 &  36.0  \\
3 & $A_7 \sim A_8$          & 28.0 &  36.0  \\
4 & $A_9 \sim A_{10}$     	& 26.0 &  36.0  \\
5 & $A_{11} \sim A_{14}$	& 22.0 &  36.0\\
6 & $A_{15}\sim A_{18}$	& 20.0 &  36.0\\
7 & $A_{19} \sim A_{22}$	& 18.0 &  36.0\\ \hline
\label{T22bar1}
\end{tabular}
\end{table}
\begin{table}
\caption{Loading conditions for the 22-bar space truss}
\centering
\begin{tabular}{l*{4}{c}r}
Node & Loading Condition 1(psi) & Loading Condition 2(psi) & Loading Condition 3 (psi) \\ \hline
Direction & $P_x , P_y , P_z$ & $P_x , P_y , P_z$ & $P_x , P_y , P_z$ \\ \hline
1 & -20.0,0.0,-5.0 & -20.0,-5.0,0.0 & -20.0,0.0,35.0&\\
2 & -20.0,0.0,-5.0 & -20.0,-50.0,0.0& -20.0,0.0,0.0&\\
3 & -20.0,0.0,-30.0& -20.0,-5.0,0.0 & -20.0,0.0,0.0&\\
4 & -20.0,0.0,-30.0& -20.0,-50.0,0.0& -20.0,0.0,-35.0&\\ \hline
\label{T22bar2}
\end{tabular}
\end{table}
\begin{table}
\caption{Optimal design comparison for the 22-bar planar truss}
\centering
\resizebox{\textwidth}{!}{%
\begin{tabular}{l*{5}{c}r}
Variables($in^2$)   & Sheu and Schmit \cite{Sheu1972} & Khan and Willmert \cite{Khan1979} & Farshi and Ziazi  \cite{Farshi2010} & This work  \\ \hline
$A_1$ 			& 2.629 & 2.563 & 2.6250 & 2.6301\\
$A_2$                                & 1.162 & 1.553 & 1.2164 & 1.2289\\
$A_3$                                & 0.343 & 0.281 & 0.3466 & 0.3550\\
$A_4$     		          & 0.423 & 0.512 & 0.4161 & 0.4153\\
$A_5$			          & 2.782 & 2.626 & 2.7732 & 2.7332\\
$A_6$			          & 2.173 & 2.131 & 2.0870 & 2.0688\\
$A_7$			          & 1.952 & 2.213 & 2.0314 & 2.0371\\
Weight($lb$)		          & 1024.8 & 1034.74 & 1023.9 & 1019.43\\ \hline
Note: 1 $in^2 = 6.425 cm^2$ $1 lb = 4.448 N$
\label{T22bar3}
\end{tabular}}
\end{table}
\subsection{The 25-bar space truss}
\begin{figure}[thb]
\begin{center}
  \centering
  \includegraphics[clip, width = 120mm , height = 75mm]{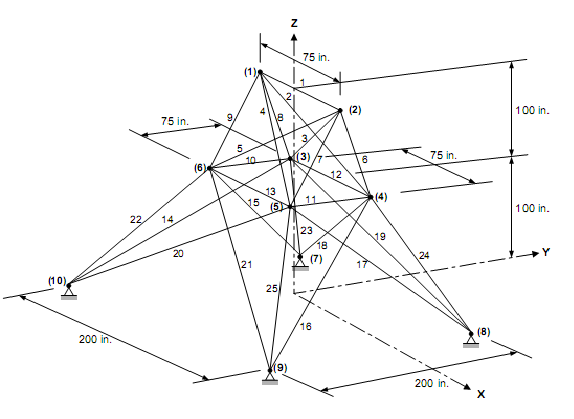} 
\caption{25-bar space truss}
\label{fig:25bar}
\end{center}
\end{figure}

\begin{figure}[thb]
\begin{center}
  \centering
  \includegraphics[clip, width = 120mm , height = 75mm]{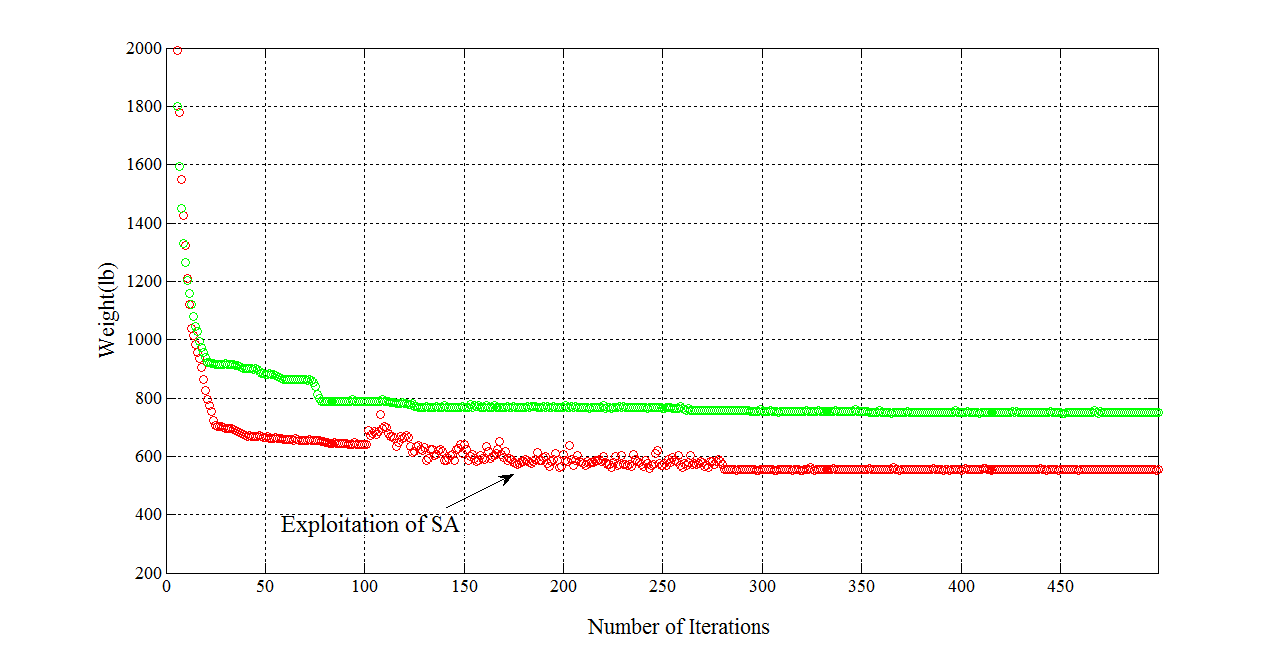} 
\caption{The 25-bar convergence history}
\label{fig:25barConv}
\end{center}
\end{figure}

The 25-bar tower space truss of Fig.\ref{fig:25bar} has been analyzed by many researchers, including Schmit and Farshi \cite{L.A.Schmit1974}, Schmit and Miura \cite{L.A.Schmit1976}, Venkayya \cite{Venkayya1971},  Adeli and Kamal \cite{Adeli1986}, Saka \cite{Saka1990}, Lamberti \cite{Lamberti2008}, Farshi and Ziazi \cite{Farshi2010}. In these studies, the material density is $0.1 lb/(in^3 )$ and modulus of elasticity is $10000 ksi$. This space truss is subjected to the two loading conditions shown in Table \ref{TT25bar}. 
The cross section grouping is as follows: $(1) A_1  (2) A_2 \sim A_5  (3) A_6 \sim A_9  (4) A_{10} \sim A_{11} (5) A_{12} \sim A_{13}  (6) A_{14} \sim A_{17}  (7) A_{18} \sim A_{2}1  (8) A_{22} \sim A_{25}$
\begin{table}
\caption{Loading conditions for the 25-bar space truss}
\centering
\begin{tabular}{l*{4}{c}r}
Node & Loading Condition 1(kips) & Loading Condition 2(kips) \\ \hline
 & $P_x , P_y , P_z$ & $P_x , P_y , P_z$ \\ \hline
1 & 0.0,20.0,-5.0 & 1.0,10.0,-5.0 \\
2 & 0.0,-20.0,-5.0 & 0,10.0,-5.0 \\
3 & 0.0,0.0,0.0 & 0.5,0.0,0.0 \\
6 & 0.0,0.0,0.0 & 0.5,0.0,0.0\\ \hline
\label{TT25bar}
\end{tabular}
\end{table}
\begin{table}
\caption{Member stress limitations for the 25-bar space truss}
\centering
\begin{tabular}{l*{3}{c}r}
Variables($in^2$)   & Compressive Stress Limitation (ksi) &   Tensile stress limitation(ksi) \\ \hline
$A_1$ 				& 35.092 & 40.0\\
$A_2 \sim A_5$       & 11.590 & 40.0\\
$A_6 \sim A_9$       & 17.305 & 40.0\\
$A_{10} \sim A_{11}$   & 35.092 & 40.0\\
$A_{12} \sim A_{13}$	& 35.092 & 40.0\\
$A_{14}\sim A_{17}$		& 6.759 & 40.0\\
$A_{18} \sim A_{21}$	& 6.959 & 40.0\\
$A_{22} \sim A_{25}$	& 11.082 & 40.0\\ \hline
\label{T25bar2}
\end{tabular}
\end{table}
The truss members were subjected to the compression and tensile stress limitations shown in Table \ref{T25bar2} In addition, maximum displacement limitations of $\pm 0.35$ in were imposed on every node in every direction. The minimum cross-sectional area of all members was $0.01 in^2$. In Fig.\ref{fig:25barConv} the convergence histories of a simple GA and H-SAGA are shown for the 25-bar truss. It could be seen that H-SAGA converges much faster and with fewer function evaluations and leads to better solution vectors compared to a simple GA. The neighbourhood search takes effect in as early as $20th$ generation, where it helps the algorithm detour from the course of the GA and find noticeably lower function values. Table \ref{T25bar3} compares optimization results with those reported in literature.
\begin{table}
\caption{Optimal design comparison for the 25-bar planar truss}
\centering
\resizebox{\textwidth}{!}{ %
\begin{tabular}{l*{9}{c}r}
Variables($in^2$)   & Schmit and Farshi   \cite{L.A.Schmit1974} & Schmit and Miura \cite{L.A.Schmit1976} & Venkayya  \cite{Venkayya1971} & Adeli and Kamal \cite{Adeli1986} & Saka \cite{Saka1990} & Lamberti \cite{Lamberti2008} & Farshi and Ziazi \cite{Farshi2010} & This work  \\ \hline
$A_1$ 			& 0.01 & 0.01 & 0.028 & 0.010 & 0.010 & 0.0100 & 0.0100 & 0.0100\\
$A_2$                                & 1.964 & 1.985 & 1.964 & 1.986 & 2.085 & 1.9870 & 1.9981 & 1.9864\\
$A_3$                                & 3.033 & 2.996 & 3.081 & 2.961 & 2.988 & 2.9935 & 2.9828 & 2.9975\\
$A_4$     		          & 0.01 & 0.01 & 0.01 & 0.010 & 0.010 & 0.0100 & 0.0100 & 0.0100\\
$A_5$			          & 0.01 & 0.01 & 0.01 & 0.010 & 0.010 & 0.0100 & 0.0100 & 0.0100\\\
$A_6$			         &  0.67 & 0.684 & 0.693 & 0.806 & 0.696 & 0.6840 & 0.6837 & 0.6806\\
$A_7$			         & 1.68 & 1.667 & 1.678 & 1.680 & 1.670 & 1.6769 & 1.6750 & 1.6733\\
$A_8$			         & 2.67 & 2.662 & 2.627 & 2.530 & 2.592	 & 2.6621 & 2.6668 & 2.6638\\ \hline
Weight($lb$)		         & 545.22 & 545.17 & 545.49 & 545.66 & 545.23 & 545.16 & 545.37 & 544.88\\ \hline
Note: 1 $in^2 = 6.425 cm^2$ $1 lb = 4.448 N$
\label{T25bar3}
\end{tabular}}
\end{table}
\subsection{The 72-bar space truss}
\begin{figure}[thb]
\begin{center}
  \centering
  \includegraphics[clip, width = 75mm , height = 75mm]{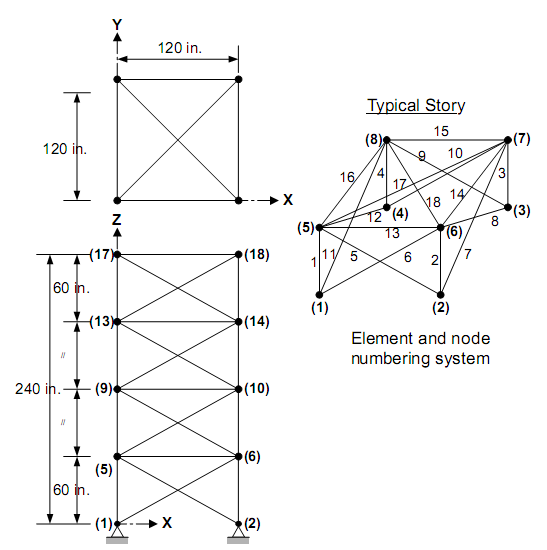} 
\caption{72-bar space truss}
\label{fig:72bar}
\end{center}
\end{figure}

The spatial 72-bar truss shown in Fig. \ref{fig:72bar} is a test case analyzed by many researchers (Schmit and Farshi \cite{L.A.Schmit1974}, Schmit and Miura \cite{L.A.Schmit1976}, Ven-kayya \cite{Venkayya1971}, Arora and Hauge \cite{J.SArora1976}, Chao et al. \cite{Chao1984}, Sedaghati \cite{Sedaghati2005}, Farshi and Ziazi \cite{Farshi2010}). Material density is $0.1 lb/in^3 $ while the modulus of elasticity is $10000 ksi$. The structure is subjected to two independent loading conditions.
\begin{enumerate}
\item Loading of $P_x=5.0 kips$, $P_y=5.0 kips$ and $P_z=-5.0 kips$  imposed on node 17; 
\item Loading of $P_z=-5.0 kips$ is applied to nodes $17, 18, 19$ and $20$;
\end{enumerate}
The loading condition divides the 72 member truss to these sub groupings:
$(1) A_1 \sim A_4 (2) A_5 \sim A_{12} (3)  A_{13} \sim A_{16}    (4) A_{17} \sim A_{18}    (5) A_{19} \sim A_{22}      (6) A_{23} \sim A_{30}  (7) A_{31} \sim A_{34} (8) A_{35} \sim A_{36}     (9) A_{37} \sim A_{40}    (10) A_{41} \sim A_{48}  (11) A_{49} \sim A_{52}  (12) A_{53} \sim A_{54}   (13) A_{55} \sim A_{58} (14) A_{59} \sim A_{66}  (15) A_{67} \sim A_{70}  (16) A_{71} \sim A_{72}$.

The members were subjected to stress limitations of $\pm 25 ksi$ and the maximum displacement of uppermost nodes was not allowed to exceed $\sim 0.25 in$ in the $x$ and $y$ directions. In this case, the minimum cross-sectional area of all members was $0.1 in^2$.\\
Optimization results are summarized in Table \ref{T72bar}. It can be seen that the proposed algorithm outperforms other works presented in the literature.
\begin{table}
\caption{Optimal design comparison for the 72-bar space truss}
\centering
\resizebox{\textwidth}{!}{%
\begin{tabular}{l*{9}{c}r}
Variables($in^2$)   & Schmit and Farshi   \cite{L.A.Schmit1974} & Schmit and Miura \cite{L.A.Schmit1976} & Venkayya  \cite{Venkayya1971} & Arora and Haug \cite{J.SArora1976} & Chao et al. \cite{Chao1984} & Sedaghati \cite{Sedaghati2005} & Farshi and Ziazi \cite{Farshi2010} & This work  \\ \hline
$A_1$          & 0.158 & 0.157 & 0.161 & 0.1564 & 0.157 & 0.1565 & 0.1565 & 0.1563\\
$A_2$          & 0.594 & 0.546 & 0.557 & 0.5464 & 0.549 & 0.5456 & 0.5457 & 0.5462\\
$A_3$          & 0.341 & 0.411 & 0.377 & 0.4110 & 0.406 & 0.4104 & 0.4106 & 0.4096\\
$A_4$          & 0.608 & 0.570 & 0.506 & 0.5712 & 0.555 & 0.5697 & 0.5697 & 0.5696\\
$A_5$	         & 0.264 & 0.523 & 0.611 & 0.5263 & 0.513 & 0.5237 & 0.5237 & 0.5239\\
$A_6$	         & 0.548 & 0.517 & 0.532 & 0.5178 & 0.529 & 0.5171 & 0.5171 & 0.5159 \\
$A_7$	         & 0.100 & 0.100 & 0.100 & 0.1000 & 0.100 & 0.1000 & 0.1000 & 0.1002\\
$A_8$	         & 0.151 & 0.100 & 0.100 & 0.1000 & 0.100 & 0.1000 & 0.1000 & 0.1006\\
$A_9$         & 1.107 & 1.267 & 1.246 & 1.2702 & 1.252 & 1.2684 & 1.2685 & 1.2691\\
$A_{10}$    & 0.579 & 0.512 & 0.524 & 0.5124 & 0.524 & 0.5117 & 0.5118 & 0.5101\\
$A_{11}$    & 0.100 & 0.100 & 0.100 & 0.1000 & 0.100 & 0.1000 & 0.1000 & 0.1000\\
$A_{12}$    & 0.100 & 0.100 & 0.100 & 0.1000 & 0.100 & 0.1000 & 0.1000 & 0.1012\\
$A_{13}$    & 2.078 & 1.885 & 1.818 & 1.8656 & 1.832 & 1.8862 & 1.8864 & 1.8861\\
$A_{14}$    & 0.503 & 0.513 & 0.524 & 0.5131 & 0.512 & 0.5123 & 0.5122 & 0.5129\\
$A_{15}$    & 0.100 & 0.100 & 0.100 & 0.1000 & 0.100 & 0.1000 & 0.1000 & 0.1000\\
$A_{16}$    & 0.100 & 0.100 & 0.100 & 0.1000 & 0.100 & 0.1000 & 0.1000 & 0.1009\\ \hline
Weight($lb$)	 & 388.63 & 379.64 & 381.2 & 379.62 & 379.62 & 379.62 & 379.65 & 379.56\\ \hline
Note: 1 $in^2 = 6.425 cm^2$ $1 lb = 4.448 N$
\label{T72bar}
\end{tabular}}
\end{table}
\subsubsection{The 200 bar truss structure with three independent loading conditions \label{200bar-1}}
\begin{figure}[thb]
\begin{center}
  \centering
  \includegraphics[clip, width = 90mm , height = 150mm]{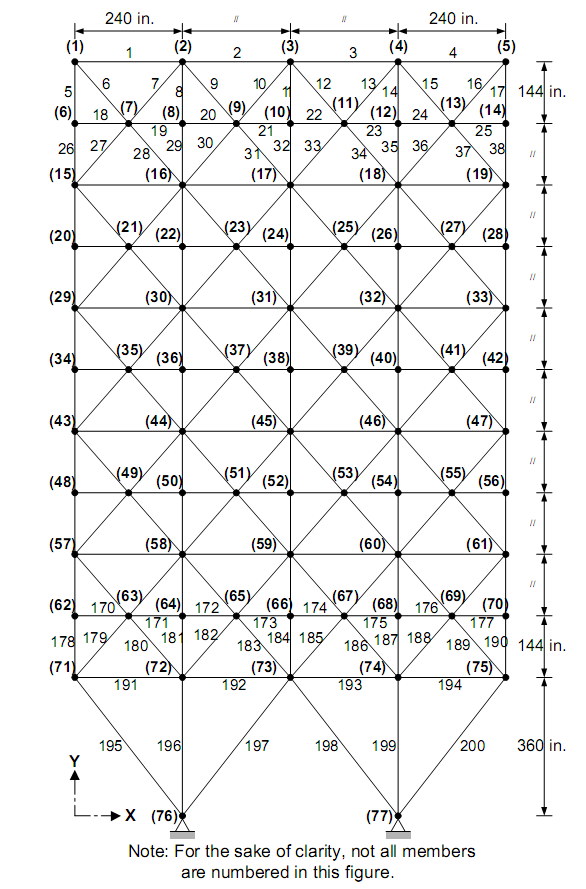} 
\caption{200-bar space truss}
\label{fig:200bar}
\end{center}
\end{figure}

\begin{figure}[thb]
\begin{center}
  \centering
  \includegraphics[clip, width = 150mm , height = 75mm]{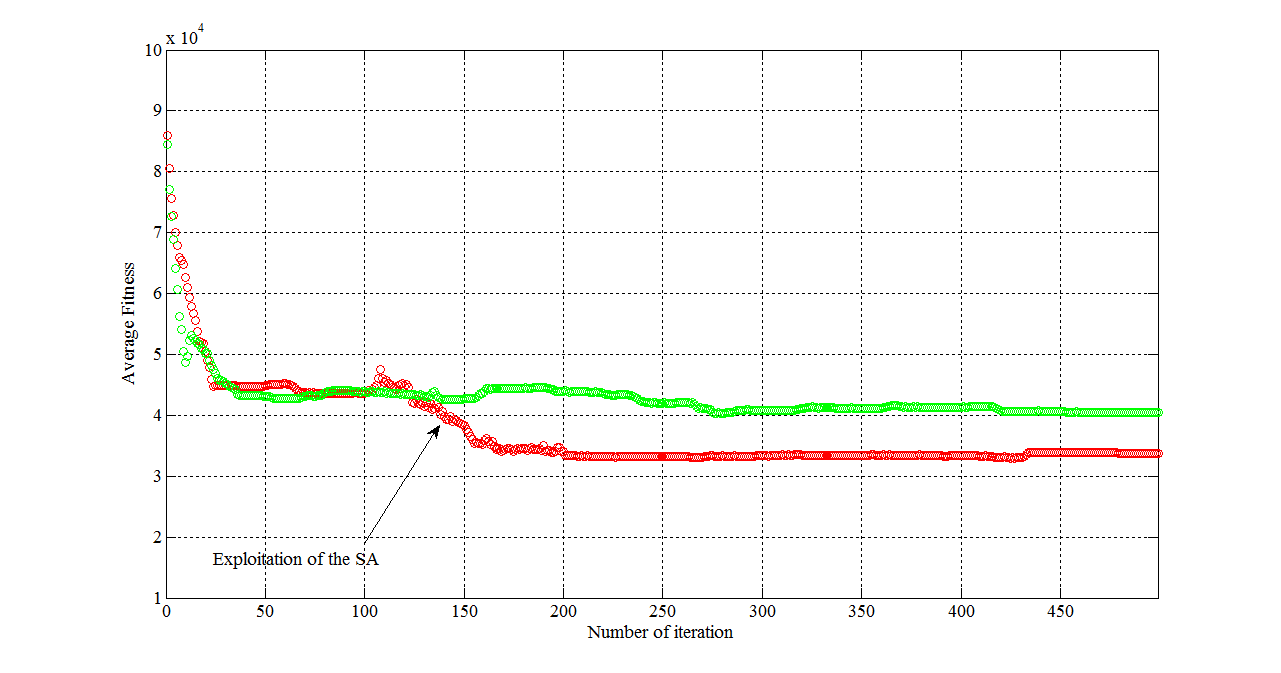} 
\caption{The 200-bar convergence history}
\label{fig:200bar2}
\end{center}
\end{figure}

The planar 200-bar truss, shown in Fig.\ref{fig:200bar}, sizing optimization problem was solved by Lamberti \cite{Lamberti2008} and Farshi and Ziazi  \cite{Farshi2010}. The steel members have a density of  $0.283 lb/in^3$ and modulus of elasticity $30000 ksi$, respectively. There are only constraints on member stresses that must be lower than $10000 psi$ (Same limit in tension and compression). The structure is subjected to three independent loading condition:\\ 
$(a)1000 lbf$ in positive $x-direction$ at nodes 1, 6, 15, 20, 29, 43, 48, 57, 62 and 71;\\
$(b)10000 lbf$  in negative $y-direction$ at nodes 1, 2, 3, 4, 5, 6, 8, 10, 12, 14, 15, 16, 17, 18, 19, 20, 22, 24, 26, 28, 29, 30, 31, 32, 33, 34, 36, 38, 40, 42, 43, 44, 45, 46, 47, 48, 50, 52, 54, 56, 58, 59, 60, 61, 62, 64, 66, 68, 70, 71, 72, 73, 74 and 75;\\
$(c)$ The loading conditions (a) and (b) acting together.\\
Groupings were done in order to have 200 independent variables. Table \ref{T200bar1} compares optimization results with literature.
In Fig.\ref{fig:200bar2}  the convergence histories of a simple GA and H-SAGA are shown for the 200-bar truss. Again, it could be seen that H-SAGA converges much faster and with fewer function evaluations and leads to better solution vectors. This time the detour begins in about $120^{th}$ generation.
\clearpage
\begin{table}
\caption{Optimal design comparison for the 200-bar planar truss}
\centering
\resizebox{\textwidth}{!}{%
\begin{tabular}{l*{5}{c}r}
Variables($in^2$)   & Members  & Lamberti \cite{Lamberti2008} & Farshi and Ziazi  \cite{Farshi2010} & This work  \\ \hline
$A_1$ 		& 1,2,3,4 & 0.1467 & 0.147 & 0.1457\\
$A_2$                     & 5,8,11,14,17 & 0.94 & 0.945 & 0.9405\\
$A_3$                     & 19,20,21,22,23,24 & 0.1 & 0.1	 & 0.1004\\
$A_4$     	          & 18,25,56,63,94,101,132,139,170,177 & 0.1 & 0.1 & 0.1 \\
$A_5$		          & 26,29,32,35,38 & 1.94 & 1.9451 & 1.9397 \\
$A_6$		          & 6,7,9,10,12,13,15,16,27,28,30,31,33,34,36,37 & 0.2962 & 0.2969 & 0.2958\\
$A_7$		          & 39,40,41,42 & 0.1 & 0.1 & 0.101\\
$A_8$		          & 43,46,49,52,55 & 3.104 & 3.1062 & 3.1032\\
$A_9$ 	          & 57,58,59,60,61,62 & 0.1 & 0.1 & 0.1012\\
$A_{10}$               & 64,67,70,73,76 & 4.104 & 4.1052 & 4.1084\\
$A_{11}$               & 44,45,47,48,50,51,53,54,65,66,68,69,71,72,74,75 & 0.4034 & 0.4039 & 0.4042\\
$A_{12}$     	         & 77,78,79,80 & 0.1922 & 0.1934	 & 0.1872\\
$A_{13}$	         & 81,84,87,90,93 & 5.4282 & 5.4289 & 5.4329\\
$A_{14}$	         & 95,96,97,98,99,100 & 0.1 & 0.1 & 0.1018\\
$A_{15}$              & 102,105,108,111,114 & 6.4282 & 6.4289 & 6.4244\\
$A_{16}$	         & 82,83,85,86,88,89,91,92,103,104,106,107,109,110,112,113 & 0.5738 & 0.5745 & 0.5723\\
$A_{17}$	         & 115,116,117,118 & 0.1325 & 0.1339 & 0.1327\\
$A_{18}$              & 119,122,125,128,131 & 7.9726 & 7.9737 & 7.9708\\
$A_{19}$              & 133,134,135,136,137,138 & 0.1 & 0.1 & 0.1007\\
$A_{20}$     	        & 140,143,146,149,152 & 8.9726 & 8.9737 & 8.9735\\
$A_{21}$	        & 120,121,123,124,126,127,129,130,141,142,144,145,147,148,150,151 & 0.7048 & 0.7053 & 0.7048\\
$A_{22}$	        & 153,154,155,156 & 0.4202 & 0.4215 & 0.4192\\
$A_{23}$	        & 157,160,163,166,169 & 10.8666 & 10.8675 & 10.8671\\
$A_{24}$	        & 171,172,173,174,175,176 & 0.1	 & 0.1 & 0.1002\\
$A_{25}$             & 178,181,184,187,190 & 11.8666 & 11.8674 & 11.8649\\
$A_{26}$             & 158,159,161,162,164,165,167,168,179,180,182,183,185,186,188,189 & 1.0344 & 1.0349 & 1.0333\\
$A_{27}$     	       & 191,192,193,194 & 6.6838 & 6.6849 & 6.6852 \\
$A_{28}$	       & 195,197,198,200 & 10.8083 & 10.8101 & 10.8036\\
$A_{29}$	       & 196,199 & 13.8339 & 13.8379 & 13.8328\\ \hline
Weight($lb$)	       & 25446.76 & 25456.57 & 25443.11\\ \hline
Note: 1 $in^2 = 6.425 cm^2$ , $1 lb = 4.448 N$
\label{T200bar1}
\end{tabular}}
\end{table}
\subsubsection{The 200 bar truss structure with five independent loading conditions}
The planar 200-bar truss structure shown in Fig.\ref{fig:200bar} can be optimized also with 200 design variables by assigning a sizing variable to the cross-sectional of each element. The structure is designed to carry five independent loading conditions. Besides the three load cases listed in Section \ref{200bar-1}, the structure is also loaded by $(d)1000lbf$ acting in the negative $x-direction$ at node points 5, 14, 19, 28, 33, 42, 47, 56, 61, 70 and 75;\\
$(e)$ Loading conditions $(b)$ and $(d)$ described in section \ref{200bar-1} acting together.\\
The optimization problem includes 3500 non-linear constraints on nodal displacements and member stresses. The displacements of all free nodes in both directions $x$ and $y$ must be less than $\pm 0.5in$. The allowable stress (the same in tension and compression) is $30000psi$. The lower bound of cross-sectional areas is $0.1 in^2$. The latest and the best optimal result reported by Lamberti and Pappalettere\cite{Lamberti2003} is $28781 lb$ but H-SAGA outperforms it by obtaining $28661 lb$ for the weight of this structure.

\subsection{The 112-bar dome truss}
\begin{figure}[thb]
\begin{center}
  \centering
  \includegraphics[clip, width = 100mm , height = 100mm]{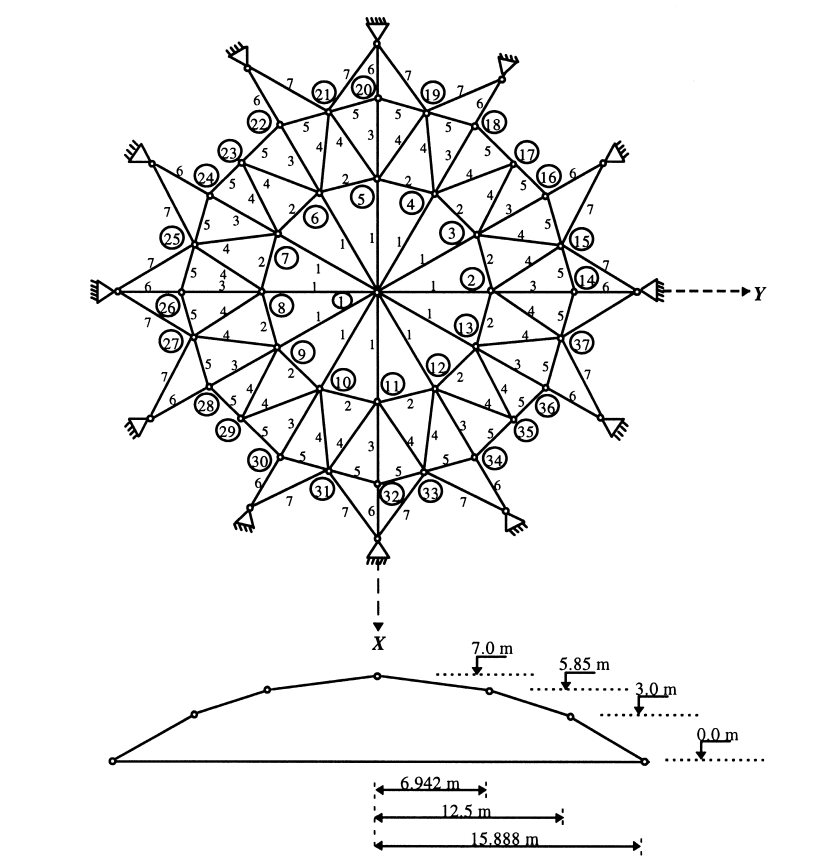} 
\caption{The 112-bar dome truss}
\label{fig:112bar}
\end{center}
\end{figure}

Fig.\ref{fig:112bar} shows the 112-bar steel dome that has been previously discussed in Saka \cite{Saka1990} and Erbatur and Hasancebi  \cite{Erbatur2000}. The optimization problem is modified as follows. First, 112 members of the dome are collected in seven distinct groups, while two groups were used in Saka\cite{Saka1990}. Next, in place of $AISC$ specification, the allowable compressive stress for each member is computed according to the Turkish specification. Despite the fact that the Turkish specification leads to safer values in computing the allowable compressive stresses, the distinction is not of much consequence. Loading conditions and optimization results are respectively presented in Table \ref{T122barl} and \ref{T112bar2}. Complete detailed list of design data is presented as follow:\\
Displacement constraints: $\delta_j \leq 20 mm$ in $z$ direction, $j=1,17,23$; \\ 
Stress constraints: $\sigma_{tention_i} \leq 34.8 ksi$,  $\sigma_{compression_i} \leq TURKISH spec, i=1,...,112$;\\
\begin{table}
\caption{Loading condition for the 112-bar dome}
\centering
\begin{tabular}{l*{4}{c}r}
Joint's number   & $x (kips)$ & $y (kips)$ & $z (kips)$ \\ \hline
1 &  0 		 & 5 & 12  \\
17,23,29,35 &   0   & 9 &  9  \\
16,18,22,24, 28,30,31,32    & 8 &  7 &  7  \\
others & 8 &  7 & 9  \\  \hline
\label{T122barl}
\end{tabular}
\end{table}
\begin{table}
\caption{Optimal design comparison for the 112- bar dome truss}
\centering
\begin{tabular}{l*{3}{c}r}
Variables($in^2$)   & Erbatur  \cite{Erbatur2000} & This work  \\ \hline
$A_1$ 	& 1.095 & 1.059\\
$A_2$           & 0.863 & 0.888\\
$A_3$           & 1.033 & 0.929\\
$A_4$     	& 1.095 & 0.948\\
$A_5$		& 1.095 & 1.024\\
$A_6$		& 0.81 & 0.9\\
$A_7$		& 1.735 & 1.612\\ \hline
Weight($lb$)	 & 7657.778 & 7015.053\\ \hline
Note: 1 $in^2 = 6.425 cm^2$ $1 lb = 4.448 N$
\label{T112bar2}
\end{tabular}
\end{table}
\subsection{The 132-bar geodesic dome}
\begin{figure}[thb]
\begin{center}
  \centering
  \includegraphics[clip, width = 100mm , height = 75mm]{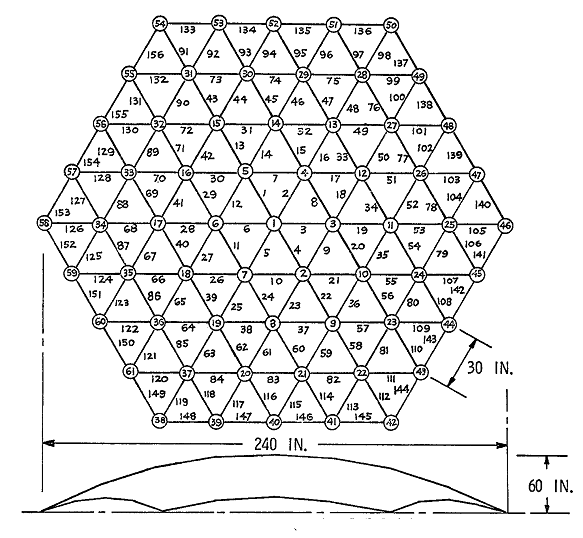} 
\caption{The 132-bar dome truss}
\label{fig:132bar}
\end{center}
\end{figure}

Fig. \ref{fig:132bar} shows a geodesic truss dome consisting of 132 bar elements that has been analyzed by Farshi and Ziazi \cite{Farshi2010}. Stress limitations for all members is $25 ksi$. Cross sectional areas must be greater than $0.1in^2$. The areas are linked together to form 36 independent design variables (Table \ref{T132bar2}). Bottom nodes are fixed to ground in all coordinate directions while free nodes cannot displace by more than $0.1 in$. The structure is subject to four independent loading conditions (Table \ref{T132bar1}). Optimization results are summarized in Table \ref{T132bar2}.

\begin{table}
\caption{Loading conditions for the 132-bar geodesic dome}
\centering
\begin{tabular}{l*{2}{c}r}
Loaded joints  & 1000 lb downward acts on each node\\ \hline
1 & $A_1\sim A_4$  \\
2 & $A_5 \sim A_6$       \\
3 & $A_7 \sim A_8$        \\
4 & $A_9 \sim A_{10}$    \\ \hline
\label{T132bar1}
\end{tabular}
\end{table}

\begin{table}
\caption{Optimal design comparison for the 132-bar geodesic dome}
\centering
\resizebox{\textwidth}{!}{%
\begin{tabular}{l*{8}{c}r}
Variables($in^2$)   & Members   & Farshi and Ziazi \cite{Farshi2010} & This work   & Variables  & Members & Farshi and Ziazi \cite{Farshi2010} & This work  \\ \hline
$A_1$ 		& 3,6 & 0.9876 & 0.9915&  $A_{19}$  & 44,47,59,62 & 0.3910 & 0.3834\\
$A_2$                     & 1,2,4,5 & 0.9902 & 0.9834 & $A_{20}$  & 45,46,60,61 & 0.4597 & 0.455\\
$A_3$                     & 8,9,11,12 & 0.9041 & 0.899 & $A_{21}$ & 78,79,87,88 & 0.2128 & 0.212\\
$A_4$     	          & 7,10	& 0.8509 & 0.8495 &  $A_{22}$  & 77,80,86,89 & 0.1721 & 0.1713\\
$A_5$		          & 19,28 & 0.3703 & 0.3692 & $A_{23}$ & 76,81,85,90 & 0.2264 & 0.2275\\
$A_6$		          & 18,20,27,29 & 0.4920 & 0.4852 & $A_{24}$ & 73,75,82,84 & 0.1657 & 0.1624\\
$A_7$		          & 17,21,26,30 & 0.4145 & 0.4101 & $A_{25}$ & 74,83 & 0.1000 & 0.101\\
$A_8$		          & 13,16,22,25 & 0.5519 & 0.5452 & $A_{26}$ & 105,126 & 0.1000 & 0.1041\\
$A_9$ 	          & 14,15,23,24 & 0.5110 & 0.5052 & $A_{27}$ &104,106,125,127 & 0.3229 & 0.3223\\
$A_{10}$                & 34,35,40,41 & 0.4278 & 0.4225 & $A_{28}$ & 103,107,124,128 & 0.3228 & 0.3167\\
$A_{11}$               & 33,36,39,42 & 0.4565 & 0.459 & $A_{29}$ & 102,108,123,129 & 0.4175 & 0.4144\\
$A_{12}$     	          & 31,32,37,38 & 0.3808 & 0.3722 & $A_{30}$ & 101,109,122,130 & 0.2587 & 0.256\\
$A_{13}$	          & 53,68 & 0.3115 & 0.3087 & $A_{31}$ & 100,110,121,131 & 0.4973 & 0.4911\\
$A_{14}$	          & 52,54,67,69 & 0.3826 & 0.3772 & $A_{32}$ & 99,111,120,132 & 0.2660 & 0.267\\
$A_{15}$	          & 51,55,66,70 & 0.3804 & 0.375 & $A_{33}$ & 91,98,112,119 & 0.1000 & 0.1001\\
$A_{16}$	         & 50,56,65,71 & 0.4412 & 0.4385 & $A_{34}$ & 92,97,113,118 & 0.3394 & 0.338\\
$A_{17}$	         & 49,57,64,72 & 0.3359 & 0.3337 & $A_{35}$ & 93,96,114,117 & 0.3425 & 0.3434\\
$A_{18}$	        & 43,48,58,63	& 0.5058 & 0.5083 & $A_{36}$ & 94,95,115,116 & 0.2986 & 0.2967\\ 
$ $                         & & & & & Weight ($lb$) & 172.74 & 171.52\\ \hline
Note: 1 $in^2 = 6.425 cm^2$ $1 lb = 4.448 N$
\label{T132bar2}
\end{tabular}}
\end{table}
\clearpage
\subsection{The 26-story tower truss}
\begin{figure}[thb]
\begin{center}
  \centering
  \includegraphics[clip, width = 75mm , height = 100mm]{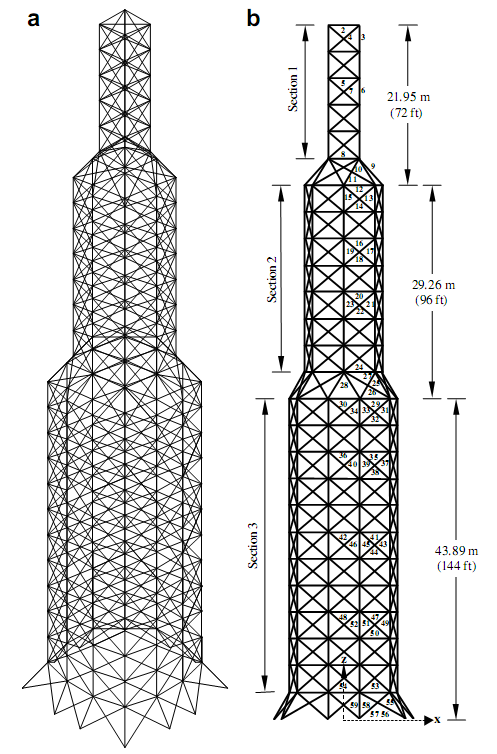} 
\caption{The 26-Story, 942-bar space truss tower: (a) 3D view, and (b) Side view }
\label{fig:26sbar}
\end{center}
\end{figure}

\begin{figure}[thb]
\begin{center}
  \centering
  \includegraphics[clip, width = 75mm , height = 75mm]{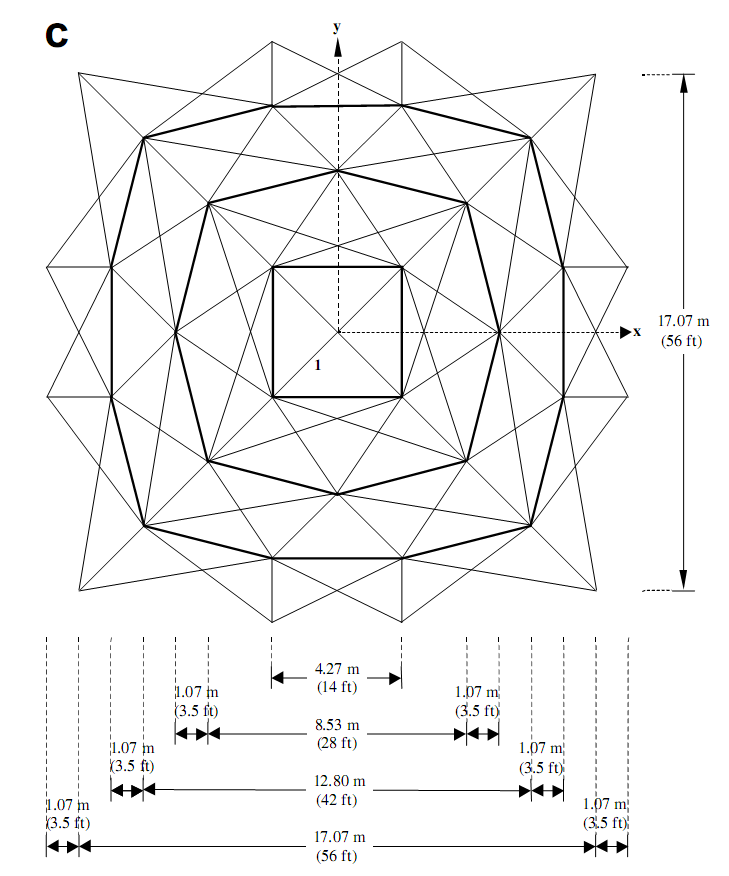} 
\caption{The 26-Story, 942-bar space truss tower: (c) Top view}
\label{fig:26sbar2}
\end{center}
\end{figure}

The 26-story-tower space truss containing 942 elements and 244 nodes is considered in this section. 59 design variables are used to represent the cross-sectional areas of 59 element groups in this structure because of structural symmetry. Figures \ref{fig:26sbar} and \ref{fig:26sbar2} show the geometry and the 59 element groups. The material density is $0.1 lb/in^3 $ and the modulus of elasticity is $10000 ksi$. The members are subject to the stress limit of $\pm 25 ksi$ ($\pm 172.375 MPa$) and the four nodes of the top level in the $x$, $y$ and $z$ directions are subject to the displacement limits of $\pm 15.0 in$ ($\pm 38.10 cm$) (about $1/250$ of the total height of the tower). The allowable cross-sectional areas in this example are selected from $0.1$ to $20.0 in^2$ ($0.6452 ~ 129.03 cm^2$). The loading on the structure consists of:\\
(1) The vertical load at each node in the first section is equal to $3 kips$; \\
(2) The vertical load at each node in the second section is equal to $6 kips$;\\
(3) The vertical load at each node in the third section is equal to $9 kips$; \\
(4) The horizontal load at each node on the right side in the $x$ direction is equal to $1 kips$; \\
(5) The horizontal load at each node on the left side in the $x$ direction is equal to $1.5 kips$; \\
(6) The horizontal load at each node on the front side in the $y$ direction is equal to $1 kips$;\\
(7) The horizontal load at each node on the back side in the $y$ direction is equal to $1 kips$;\\
Table \ref{T26sbar} compares optimization results with literature.
\begin{table}
\caption{Optimal design comparison for the 26-story tower truss }
\centering
\begin{tabular}{l*{4}{c}r}
Cross Section($in^2$)   &  Erbatur and Hasancebi   \cite{Erbatur2000} & Hasancebi   \cite{Hasancebi2008} & This work  \\ \hline
$A_1$          & 1 & 1.02 & 3.3409\\
$A_2$          & 1 & 1.037 & 1.0226\\
$A_3$          & 3 & 2.943 & 5.7605\\
$A_4$          & 1 & 1.92 & 2.4798\\
$A_5$	          &1 & 1.025 & 1.0000\\
$A_6$	          &17 & 14.961 & 14.3539\\
$A_7$	          & 3 & 3.074 & 2.8752\\
$A_8$	          & 7 & 6.78 & 11.7044\\
$A_9$          & 20 & 18.58 & 14.8708\\
$A_{10}$     & 1 & 2.415 & 3.6599\\
$A_{11}$     & 8 & 6.584 & 5.1913\\
$A_{12}$     & 7 & 6.291 & 5.5767\\
$A_{13}$     & 19 & 15.383 & 14.1554\\
$A_{14}$     & 2 & 2.1 & 2.1912\\
$A_{15}$     & 5 & 6.021 & 2.9070\\
$A_{16}$     & 1 & 1.022 & 1.0000\\
$A_{17}$     & 22 & 23.099 & 18.1769\\
$A_{18}$     & 3 & 2.889 & 2.5274\\
$A_{19}$     & 9 & 7.96 & 12.6091\\
$A_{20}$     & 1 & 1.008 & 1.0361\\
$A_{21}$     & 34 & 28.548 &  31.1194\\
$A_{22}$     & 3 & 3.349 & 2.8803\\
$A_{23}$     & 19 & 16.144 & 17.0459\\
$A_{24}$     & 27 & 24.822 & 18.3234\\
$A_{25}$     & 42 & 38.401 & 38.7810\\
$A_{26}$     & 1 & 3.787 & 2.6226\\
$A_{27}$     & 12 & 2.32 & 9.2714\\
$A_{28}$     & 16 & 17.036 & 13.0850\\
$A_{29}$     & 19 & 14.733 & 13.5173\\
$A_{30}$     & 14 & 15.031 & 16.3403\\
$A_{31}$     & 42 & 38.597 & 37.3477\\
$A_{32}$     & 4 & 3.511 & 3.1946\\
$A_{33}$     & 4 & 2.997 & 6.3378\\
$A_{34}$     & 4 & 3.06 & 1.8829\\
$A_{35}$     & 1 & 1.086 & 1.0000\\
$A_{36}$     & 1 & 1.462 & 1.1151\\ \hline
Weight($lb$)	 & 143436 & 141241 & 140674\\ \hline
Note: 1 $in^2 = 6.425 cm^2$ $1 lb = 4.448 N$
\label{T26sbar}
\end{tabular}
\end{table}
\section{Conclusion}
A new hybrid algorithm, H-SAGA, was presented for the optimal design of planar- and space truss structures. The Genetic-based algorithm features a local-search engine that benefits from a Simulated Annealing-like stochastic search scheme. The algorithm was further ``tuned", e.g. by adopting a dynamic local search strategy, to better suit the complex nature of the problem at hand. Ten different structures, including a 112 member dome truss and a 26-story tower truss were optimally designed using the proposed algorithm and the results were compared with those reported in the literature (wherever available). The comparison showed that H-SAGA outperformed other algorithms, in some cases noticeably, both in terms of solution optimality and computational cost. In addition to its fairly high convergence rate which causes computational efficiency, the proposed algorithm is distinguished by its ability to fluently escape the traps of the local minima and to carefully avoid constraint violations when different random initial points are used.  

\clearpage

\bibliographystyle{unsrt}
\bibliography{references}

\begin{thebibliography}{10}

\bibitem{Sivanandam2008}
S.~N. Sivanandam and S.~N. Deepa.
\newblock {\em Introduction to Genetic Algorithms}.
\newblock Springer Berlin Heidelberg, New York, 2008.

\bibitem{Holland1992}
John~H. Holland.
\newblock {\em Adaptation in Natural and Artificial Systems}.
\newblock MIT press, Cambridge, 1992.

\bibitem{Bolte1996}
A.~Bolte and U.W. Thonemann.
\newblock Theory and methodology of optimizing simulated annealing schedules
  with genetic programming.
\newblock {\em European Journal of Operational Research}, 92:402 416, 1996.

\bibitem{Wang2004}
Z.G Wang, Y.S Wong, and M.~Rahman.
\newblock Optimization of multi-pass milling using genetic algorithm and
  genetic simulated annealing.
\newblock {\em The International Journal of Advanced Manufacturing Technology},
  24:727--732, 2004.

\bibitem{Wren2007}
Z.~Wren, Y.~San, and J.~F. Chen.
\newblock Hybrid simplex-improved genetic algorithm for global numerical
  optimization.
\newblock {\em Acta Automotica Sinica}, 33:91--95, 2007.

\bibitem{Rahami2011}
H.~Rahami, A.~Kaveh, M.~Aslani, and R.~Najian Asl.
\newblock A hybrid modified genetic-nelder mead simplex algorithm for
  large-scale truss optimization.
\newblock {\em Int. J. Optim. Civil Eng.}, 1:29--46, 2011.

\bibitem{Aslani2010}
M.~Aslani, R.~Najian Asl, R.~Oftadeh, and M.~Shariat Panahi.
\newblock A novel hybrid simplex-genetic algorithm for the optimum design of
  truss structures.
\newblock In {\em Proceedings of the World Congress on Engineering 2010 Vol
  II}, 2010.

\bibitem{N.Dugan2009}
S.~Erkoc N.~Dugan.
\newblock Genetic algorithm-monte carlo hybrid geometry optimization method for
  atomic clusters.
\newblock {\em Computational Materials Science}, 45:127--132, 2009.

\bibitem{Fan2006}
S.~K.~S. Fan, Y.~C. Liang, and E.~Zahara.
\newblock A genetic algorithm and a particle swarm optimizer hybridized with
  nelder�mead simplex search.
\newblock {\em Computers \& Industrial Engineering}, 50:401--425, 2006.

\bibitem{Carbas2010}
O.~Hasancebi, S.~Carbas, and M.P Saka.
\newblock Improving the performance of simulated annealing in structural
  optimization.
\newblock {\em Structural and Multidisciplinary Optimization}, 41:189--203,
  2010.

\bibitem{Kirkpatrick1983}
S.~Kirkpatrick and M.~P.~Vecchi Jr~C. D.~Gelatt.
\newblock Optimization by simulated annealing.
\newblock {\em Science Journal}, 220:671--699, 1983.

\bibitem{L.A.Schmit1974}
L.A. Schmit and B.~Farshi.
\newblock Some approximation concepts for structural synthesis.
\newblock {\em AIAA Journal}, 12:692--699, 1974.

\bibitem{L.A.Schmit1976}
L.A. Schmit and H.~Miura.
\newblock Approximation concepts for efficient structural synthesis.
\newblock {\em NASA-CR-2552}, 1:163--168, 1976.

\bibitem{Venkayya1971}
V.B Venkayya.
\newblock Design of optimum structures.
\newblock {\em Computers and Structures}, 1:265--309, 1971.

\bibitem{Sedaghati2005}
R.~Sedaghati.
\newblock Benchmark case studies in structural design optimization using the
  force method.
\newblock {\em International Journal of Solids and Structures}, 42:5848--5871,
  2005.

\bibitem{Kaveh2006}
A.~Kaveh and H.~Rahami.
\newblock Analysis, design and optimization of structures using force method
  and genetic algorithm.
\newblock {\em International Journal for Numerical Methods in Engineering},
  65:1570--1584, 2006.

\bibitem{Li2007}
L.~J. Li, Z.~B. Huang, F.~Liu, and Q.~H. Wu.
\newblock A heuristic particle swarm optimizer for optimization.
\newblock {\em Computers and Structures}, 85:340--349, 2007.

\bibitem{RizziMay1976}
P.~Rizzi.
\newblock Optimization of multi-constrained structures based on optimality
  criteria.
\newblock In {\em Proceedings of the AIAA/ASME/SAE 17th Structures Structural
  Dynamics and Materials Conference}, pages 448--462, May 1976.

\bibitem{John1987}
K.~V. John, C.~V Ramakrishnan, and KG. Sharma.
\newblock Minimum weight design of trusses using improved move limit method of
  sequential linear programming.
\newblock {\em Computers and Structures}, 27:583--591, 1987.

\bibitem{Lee2004}
K.~S. Lee and Z.~W. Geem.
\newblock A new structural optimization method based on the harmony search
  algorithm.
\newblock {\em Computers and Structures}, 82:781--798, 2004.

\bibitem{Farshi2010}
B.~Farshi and A.~Ziazi.
\newblock Sizing optimization of truss structures by method of centers and
  force formulation.
\newblock {\em International Journal of Solids and Structures}, 47:2508--2524,
  2010.

\bibitem{Venkayya1973}
V.~B. Venkayya, N.~S. Khot, and L.~Berke.
\newblock {\em Application of Optimality Criteria Approaches to Automated
  Design of Large Practical Structures}.
\newblock AGARD-CP-123, 1973.

\bibitem{H.Adeli1995}
H.~Adeli and S.~Kumar.
\newblock Distributed genetic algorithm for structural optimization.
\newblock {\em J Aerospace Eng, ASCE}, 8,3:156--163, 1995.

\bibitem{Imai1981}
K.~Imai and L.~A. Schmit.
\newblock Con?guration optimization of trusses.
\newblock {\em J Structural Division, ASCE}, 107:745--756, 1981.

\bibitem{Sheu1972}
C.~Y. Sheu and L.~A. Schmit.
\newblock Minimum weight design of elastic redundant trusses under multiple
  static load conditions.
\newblock {\em AIAA Journal}, 10:155--162, 1972.

\bibitem{Khan1979}
M.~R. Khan, K.~D. Willmert, and W.~A. Thornton.
\newblock An optimality criterion method for large-scale structures.
\newblock {\em AIAA Journal}, 17:753--761, 1979.

\bibitem{Adeli1986}
H.~Adeli and O.~Kamal.
\newblock Efficient optimization of space trusses.
\newblock {\em Computers and Structures}, 24:501--511, 1986.

\bibitem{Saka1990}
MP. Saka.
\newblock Optimum design of pin-jointed steel structures with practical
  applications.
\newblock {\em J Struct Eng. ASCE}, 116, 10:2599--2620, 1990.

\bibitem{Lamberti2008}
L.~Lamberti.
\newblock An efficient simulated annealing algorithm for design optimization of
  truss structures.
\newblock {\em Computers and Structures}, 86:1936--1953, 2008.

\bibitem{J.SArora1976}
JR.~Haug J.~S~Arora.
\newblock Ef?cient optimal design of structures by generalized steepest descent
  programming.
\newblock {\em International Journal for Numerical Methods in Engineering},
  10:747�766, 1976.

\bibitem{Chao1984}
N.~H. Chao, S.~J. Fenves, and A.~W. Westerberg.
\newblock {\em Application of reduced quadratic programming technique to
  optimal structural design}.
\newblock John Wiley, New York, 1984.

\bibitem{Lamberti2003}
Lamberti and Pappalettere.
\newblock Move limits definition in structural optimization with sequential
  linear programming, part ii: Numerical examples.
\newblock {\em Computer and Structures}, 81:215--238, 2003.

\bibitem{Erbatur2000}
F.~Erbatur and O.~Hasancebi.
\newblock Optimal design of planar and space structures with genetic
  algorithms.
\newblock {\em Computer and Structures}, 75:209--224, 2000.

\bibitem{Hasancebi2008}
O.~Hasancebi.
\newblock Adaptive evolution strategies in structural optimization: Enhancing
  their computational performance with applications to large-scale structures.
\newblock {\em Computers and Structures}, 86:119--132, 2008.

\end{thebibliography}

\end{document}